\documentclass{amsart}
\usepackage{amsmath,amsthm,amssymb}

\newtheorem{thmno}{Theorem}[section]
\newtheorem{thmlet}{Theorem}

\newtheorem{consno}[thmno]{Corollary}
\newtheorem{lmmno}[thmno]{Lemma}

\theoremstyle{definition}
\newtheorem{defn}[thmno]{Definition}
\newtheorem*{xrem}{Remark}

\numberwithin{equation}{section}

\DeclareMathOperator{\sgn}{sgn}
\newcommand{\eps}{\varepsilon}
\newcommand{\La}{\Lambda}
\newcommand{\la}{\lambda}
\newcommand\ebf{\mathbf{e}}
\newcommand\abf{\mathbf{a}}
\newcommand\xbf{\mathbf{x}}
\newcommand\kbf{\mathbf{k}}
\newcommand\nn{\mathbf{n}}
\newcommand\tbf{\mathbf{t}}
\newcommand\zero{\mathbf{0}}
\newcommand\LA{\mathbb L}

\newcommand\JJ{D}

\newcommand\aaa{\boldsymbol{\alpha}}
\begin{document}

\begin{center}
\large
\textbf{On Ces\'aro summability of Fourier series of functions from multidimensional Waterman classes}
\end{center}
\begin{center}
Alexandr Bakhvalov
\end{center}

\begin{quotation}
{\small \textsc{Abstract.}
An analogue of D. Waterman's result on the summability of the Fourier series for functions of bounded
 $\Lambda$-variation by the Ces\'aro methods of negative order 
is obtained in multidimensional case. It is proved that, unlike one-dimensional case, 
the continuity of function in the corresponding variation is essential  
for the convergence and even for the localization of the Ces\'aro  means for certain orders of these means. 

  References: 12 items.


}\end{quotation}

\section{Introduction}\label{sec1}

In 1972, D.\ Waterman\cite{v1} introduced  the $\La BV$ classes  of functions of bounded $\Lambda$-variation. 
One of these classes, the class of  functions of bounded harmonic variation,
proved to be a perfectly suitable  instrument for estimating the partial sums of  trigonometric Fourier series.
The $\La BV$ classes  were generalized
by A.~A.~Saakyan \cite{v2} for two-dimensional case and  
by A.~I.~Sablin \cite{sab2} for multidimensional case.
Certain results on convergence of trigonometric Fourier series for functions from the $\La BV$ 
classes were obtained  in these papers and in our papers \cite{owncv,ownint}. 

In 1976, Waterman\cite{watsum} applied the concept of  $\La$-variation  for studying the properties of the Ces\'aro means  
of trigonometric Fourier series. In this work,  $C\La V$ classes  of  functions continuous in $\La$-variation
were introduced. The classes $C \{ n^b \} V$  were used to obtain the convergence condition. Later on, Sablin proved that in one-dimensional case 
the classes $\{ n^b \} BV$ and $C \{ n^b \} V$ coincide, therefore, the condition of continuity in 
$\La$-variation turned out to be unimportant for this problem.

The aim of our article is to expand the results on Ces\'aro summability for multidimensional case.
In sections \ref{Saux} and \ref{Sconverge} we prove a sufficient condition for convergence of the Ces\'aro means
in terms of the Waterman classes. This condition requires continuity of a function in the corresponding variation, which was intoduced  in multidimensional case by the author\cite{owncv} and Dragoshanskii\cite{os}.

However,  the author and Dragoshanskii proved that the class $C \{ n^b \} V$ can be a proper subclass in $\{ n^b \} BV$ 
both for isotropic and unisotropic cases. Thus the following question arises again: is the continuity
in variation an essential condition for the convergence of the Ces\'aro means or not?
In section \ref{S:diverg} we answer this question. We prove that for certain orders of the means,
this condition cannot be omitted.

We are now turning to precise definitions and statements.
First, let us introduce the necessary notation.

By $C$ we denote an absolute constant; by $C(\cdot)$ we denote a value that depends
on the parameters listed in brackets 
 (they may be different in different cases). Let $\mathbb T = [-\pi,\pi]$.
For two sequences $\{ a_n \}$ and $\{ b_n \}$ we write $a_n \sim b_n$  
if there exists a finite positive limit of $\frac{a_n}{b_n}$ as $n\to\infty$.

Let $\Delta$ be an interval on $\mathbb R$.
By $\Omega(\Delta)$ we denote the set of all finite systems  $\{I_n \}$
of pairwise disjoint intervals 
such that $\overline{I_n}\subset \Delta$.
Let $E$  be a subset of $\Delta$. By $\Omega(\Delta\setminus E)$ 
denote the set of systems $\{ (a_n,b_n) \} $ from $\Omega(\Delta)$ 
such that $a_n \notin E$, $b_n \notin E$.

We say that $\Delta$ is an interval in $\mathbb R^m$ if 
$\Delta=\Delta^1 \times\dots\times \Delta^m$ where $\Delta^j$ are intervals (open, half-open or closed)
in $\mathbb R$; in other words, $\Delta$ is a parallelepiped with the edges parallel to the coordinate axes.
We also denote such an interval by $\Delta=\bigotimes\limits_{k=1}^m \Delta^k$.

Let $\{ \ebf_j \}_{j=1}^n$  be the standard basis of $\mathbb R^n$.
Consider a function $f(\xbf) = f (x^1,\dots, x^m)$ of $m$ variables.
We define the operators
$
\Delta_{\xbf,s,j} (f) = f (\xbf + s \ebf_j) - f(\xbf).
$
Let $I^k=(a^k,b^k)$  and $\abf = (a^1,\dots, a^m)$. We put 
\[
f(I) = f(I^1\times \dots\times I^{m}) = \Delta_{\abf,b^1 - a^1,1}  \circ \dots
\circ \Delta_{\abf,b^m - a^m,m} (f).
\]
The value $f(I)$ is called
\emph{the symmetric difference} of the function $f $ on $I$.

Let the set $\{1,\dots,m\}$ be divided into two non-intersecting subsets
$\gamma$ and $\xi$ containing $p$ and $m-p$ elements respectively.
We denote $|\gamma|= p$, $|\xi| = m-p$.
If $\xbf = (x^1,\dots,x^m)$,
then by $x^\gamma$ we denote the element of $\mathbb R^p$ with coordinates $x^j$, $j\in \gamma$.
For an interval $ I = \bigotimes\limits_{j=1}^m I^j$, by $I^\gamma$ denote
the interval $\bigotimes\limits_{j\in\gamma} I^j$.
By $x^\gamma+ I^\gamma$ (e.g. $\xbf+ I $)
denote the shift of the interval $I^\gamma$ on the vector $x^\gamma$.

By $f(I^\gamma, x^\xi)$ we denote the symmetric difference of  $f$ as the function of variables $x^j$, $j\in\gamma$
on $I^\gamma$ for the fixed values of $x^k$, $k\in\xi$.

\begin{defn}{}
We say that a nondecreasing sequence of positive numbers
 $\La=\{ \la_n \}$ \emph{determines a class of functions 
of bounded $\La$-variation (a Waterman class)},
if
$
\sum_{n=1}^\infty \frac{1}{\la_n} = \infty.
$
(In some works, the condition $\la_n \to \infty $ as $n\to\infty$ is added to the definition.)
Further on we consider only the sequences $\La$ of this kind.
By $\LA$ we denote the set of all such sequences $\La$.
By $\La(N)$ denote the partial sums $\sum_{k=1}^N \frac{1}{\la_k}$.
The sequence $\{ \la_k \}_{k=n+1}^\infty$ is denoted by $\La_n$.
We also write $H = \{ n \}_{n=1}^\infty$. It is clear that $H \in \LA$.
\end{defn}

\begin{defn}
Consider $\La^1,\dots, \La^m \in \LA$ and an interval $\Delta=\Delta^1 \times\dots\times \Delta^m \subset \mathbb R^m$. Then by
$(\La^1,\dots, \La^m)$-\emph{variation of a function $f(x^1,\dots,x^m) $
with respect to the variables
$x^1,\dots,x^m $
over}  $\Delta$
we mean
\[
V_{\La^1,\dots,\La^m}^{x^1,\dots,x^m} (f;\Delta) =
V_{\La^1,\dots,\La^m}^{\xbf} (f;\Delta) =
\sup_{ \{I^j_{k_j}\} \in \Omega(\Delta^j)
}
\sum_{k_1,\dots,k_m}
\frac{| f(I^1_{k^1}\times\dots\times I^m_{k_m} ) | }{\la^1_{k_1}\dots \la^m_{k_m}}.
\]
\end{defn}

Let a nonempty set $\gamma \subset \{ 1,\dots,m \}$ consist of the  elements
$j_1<\dots<j_p$ and
$\xi = \{ 1,\dots,m \}\setminus \gamma$.
By
\[
V_{\La^{\gamma}}^{x^\gamma} (f;(\Delta^\gamma,x^\xi)) =
V_{\La^{j_1},\dots,\La^{j_p}}^{x^\gamma} (f;(\Delta^\gamma,x^\xi))
\]
we denote $(\La^{j_1},\dots,\La^{j_p})$-variation of
$f$ as the function of variables $x^{j_1},\dots,x^{j_p}$
over the $p$-dimensional
parallelepiped $\Delta^\gamma=\Delta^{j_1} \times\dots\times \Delta^{j_p}$
for fixed values $x^\xi$ of other variables
(if $\xi$ is nonempty).
The parallelepipeds
$ \bigotimes\limits_{l=1}^p I^{j_l}_{k_{j_l}}$
we denote by $I^\gamma_{k^\gamma}$.
The products $ \la^{j_1}_{k_{j_1} } \dots \la^{j_p}_{k_{j_p}} $ we denote by $\la^\gamma_{k^\gamma}$.

Further, the value
\[
V_{\La^{\gamma}}^{x^\gamma} ( f;\Delta ) =
V_{\La^{j_1},\dots,\La^{j_p}}^{x^\gamma} (f;\Delta)
=\sup_{x^\xi \in \Delta^\xi}
V_{\La^{j_1},\dots,\La^{j_p}}^{x^\gamma} (f;(\Delta^\gamma,x^\xi).
\]
is called the $(\La^{j_1},\dots, \La^{j_p})$-\emph{variation of the function $f(\xbf) $
with respect to the variables
$x^\gamma$
over the interval} $\Delta$.

\begin{defn} Let $\Delta=\Delta^1 \times\dots\times \Delta^m$ be an interval in $\mathbb R^m$ and $f$ be a function on $\Delta$.
\emph{The (total) $(\La^1,\dots, \La^m)$-variation of the function $f(\xbf) $ over}
$\Delta$ is defined as
\[
V_{\La^1,\dots,\La^m} (f;\Delta)
=
\sum_{ \gamma\subset \{1,\dots,m\},\,\, \gamma\neq\varnothing} 
V_{ \La^{\gamma} }^{x^\gamma} (f;\Delta).
\]
The set of functions with finite total $(\La^1,\dots, \La^m)$-variation is called \emph{the class of functions of
bounded $(\La^1,\dots, \La^m)$-variation} and is denoted by
$(\La^1,\dots, \La^m)BV(\Delta)$.
If $\La^j= \La$ for all $j$ then be write briefly
$V_\La^{x^\gamma}$, $V_\La$, $\La BV(\Delta)$ and $\la_{k^\gamma}$.
The quantity $V_H$ is called the \emph{harmonic} variation.
\end{defn}

\begin{defn}
A function
$f$ from the class  $(\La^1,\dots, \La^m)BV(\Delta)$
is said to be \emph{continuous in $(\La^1,\dots, \La^m)$-variation on} $\Delta$,
if 
\begin{equation*}
\lim_{n\to\infty}
V_{\La^{j_1},\dots,\La^{j_{k-1}}, \La^{j_k}_n,\La^{j_{k+1}},
\dots,\La^{j_p} }^{x^\gamma} (f;\Delta)
= 0
\end{equation*}
for any nonempty $\gamma = \{ j_1,\dots,j_p \}
\subset \{1,\dots,m \}$
and any $j_k\in \gamma$. 
The set of functions continuous in $(\La^1,\dots, \La^m)$-variation on $\Delta$ is denoted by $C(\La^1,\dots,\La^m)V(\Delta)$.
\end{defn}

\begin{defn}
A point $\xbf_0$ is said to be \emph{a regular point} of a function $f(\xbf)$,
if there exist $2^m$ finite limits
\[
f(x_0^1 \pm 0, \dots, x_0^m \pm 0)
= \lim_{t^1,\dots,t^m \to +0} f(x_0^1 \pm t^1,\dots, x_0^m \pm t^m)
\]
for all combinations of signs.
For the regular point $\xbf_0$ we set 
\[
f^*(\xbf_0)= \frac{1}{2^m} \sum
f(x_0^1 \pm 0, \dots, x_0^m \pm 0).
\]
\end{defn}

In what follows,
we consider functions functions defined on $\mathbb T^m$ to be measurable and $2\pi$-periodic
 with respect to every variable (maybe, after the suitable continuation). 
Now we recall the definition of Ces\'aro means (see, e.g. \cite[Vol.1, Chap.3, \S1]{zygm}).
Let $\alpha > -1$. The values $A_n^\alpha$ are defined by the formula
\[
 \sum_{n=0}^\infty  A_n^\alpha x^n = (1-x)^{-\alpha-1}.
\]
Then the values
\[
 \sigma_n^\alpha  = \sum_{k=0}^n \frac{A^\alpha_{n-k} }{ A^\alpha_n } u_k.
\]
are called \emph{the Ces\'aro means of order $\alpha$, or $(C,\alpha)$-means,} of the series $\sum_{k=0}^\infty u_k$.
It is well known (see, e.g., \cite[Vol.1, Chap.3, (1.17)]{zygm}) that $A_n^\alpha \sim n^\alpha$.

For the trigonometric Fourier series of a function $f$ integrable on ${\mathbb T}$, 
its Ces\'aro means at a point $x$ are denoted by $ \sigma_n^\alpha  (f,x)$.
It is shown in \cite[Vol.1, Chap.3, \S5]{zygm} that 
\[
 \sigma_n^\alpha  (f,x) = \frac{1}{\pi} \int_{-\pi}^\pi f(x+t) K_n^\alpha (t)\, dt,
\]
where the functions
\[
 K_n^\alpha (t) = \frac{1}{ A_n^\alpha } \sum_{k=1}^n A_{n-k}^{\alpha-1} D_k (t) , \qquad n\in\mathbb N,
\]
are called \emph{the Ces\'aro kernels}. These functions are even, and
\begin{equation}
K_n^\alpha (t) = \frac{1}{ A_n^\alpha } \frac{ \sin \bigl[ \bigl( n + \frac{1}{2}+ \frac{\alpha}{2} \bigr) t
- \frac{\pi \alpha}{2}  \bigr] }{ \bigl( 2 \sin \frac{t}{2} \bigr)^{\alpha+1 } }
+ \frac{ 2 \theta \alpha }{ n \bigl( 2 \sin \frac{t}{2} \bigr)^2  } =K_n^{\alpha,*} (t) + R_n^\alpha (t).
\label{kernel}
 \end{equation} 
for $t\in (0,\pi]$, where $\theta = \theta(t,\alpha)$, $|\theta|<1$.
For the Ces\'aro kernel, the following properties hold:
\begin{equation}
 |K_n^\alpha (t) | \leq n+1 \leq 2n, \qquad
 |K_n^\alpha (t) | \leq B(\alpha) n^{-\alpha} |t|^{-(\alpha+1)}
\label{estker}
\end{equation}
and
\begin{equation}
	 \frac{1}{\pi } \int_{-\pi}^\pi K_n^\alpha(t)\, dt =1.
	 \label{unit}
\end{equation}

The Ces\'aro methods are defined in multidimensional case (see, e.g., \cite[part 2, chapter 2]{bookzh}) in the following way. 
Let $\aaa = (\alpha_1,\dots,\alpha_m)$ be a vector with $\alpha_j >-1$,
and let $S_\kbf (f,\xbf)$ be a rectangular partial sum of the trigonometric series of a function $f$.
The value 
\[
\sigma_{\nn}^{\aaa} (f,\xbf) = \Bigl( \prod_{j=1}^m A^{\alpha_j}_{n_j} \Bigr)^{-1} \sum_{\kbf=0}^\nn  \Bigl( \prod_{j=1}^m 
A_{n_j - k_j}^{\alpha_j -1} \Bigr) S_\kbf (f,\xbf),
\]
is called \emph{the Ces\'aro means of the order $\aaa$} of this series.
The multiple Ces\'aro means can be represented through the Ces\'aro kernels by the formula
\[
\sigma_{\nn}^{\aaa} (f,\xbf) = \frac{1}{\pi^m} \int_{\mathbb T^m} f(\xbf+\tbf) \prod_{j=1}^m K_{n_j}^{\alpha_j} (t^j) \,d\tbf.
\]
We consider the  Pringsheim convergence of Ces\'aro means, i.e. their convergence as $n_j$ tend to $+\infty$ independently.

We have mentioned the following one-dimensional result.

\begin{thmlet}[D. Waterman \cite{watsum}]\label{twatsum}
 Let $\alpha \in (-1,0)$. The Fourier series, $S[f]$, of a function $f$ of class $\{ n^{\alpha + 1 } \}BV(\mathbb T)$
is everywhere $(C,\alpha)$-bounded and is uniformly  $(C,\alpha)$-bounded on each closed interval of continuity. 
If $f$ is continuous in $\{ n^{\alpha + 1 } \}$-variation, then $S[f]$ 
is everywhere $(C,\alpha)$-summable to sum $\frac{1}{2}( f(x+0) + f(x-0)) $ and summability 
is uniform on each closed interval of continuity.
\end{thmlet}

Waterman also proved that the class $\{ n^{\alpha + 1 } \}BV(\mathbb T)$ cannot be replaced by a larger $\La BV$ class.
Later on, Sablin \cite{sab2} proved that in one-dimensional case  $\{ n^{\alpha + 1 } \}BV(\mathbb T) = C \{ n^{\alpha + 1 } \}V(\mathbb T)$ for any $-1<\alpha <0$.
Therefore, the condition $f\in  C \{ n^{\alpha + 1 } \}V(\mathbb T)$ of the second part of Theorem \ref{twatsum} can be replaced by the condition
$f\in \{ n^{\alpha + 1 } \}BV(\mathbb T)$.

In section 3, we prove the following multidimensional analogue of Theorem \ref{twatsum}.

\begin{thmno}\label{T:main}
Let $\alpha_j \in (-1,0)$ and  $\beta_j = \alpha_j+1$, $j=1,\dots,m$. 
Consider a function $f$ from the class
$( \{ n^{\beta_1 } \},\dots , \{ n^{\beta_m } \} )BV( {\mathbb T}^m)$.
Then its Fourier series is uniformly $(C,\aaa)$ bounded. 
If  $f$ is continuous in $( \{ n^{\beta_1 } \},\dots, \{ n^{\beta_m } \} )$-variation, then its Fourier series is $(C,\aaa)$ summable to $f^*( \xbf_0)$ at every regular point $\xbf_0$.
 If $f$ is continuous in a neighborhood of a compact set $K$, then summability is uniform on $K$.
\end{thmno}

For $m=2$, Dragoshanskii \cite[Theorem 1]{os} proved that 
\[ 
( \{ n^{\beta_1 } \},  \{ n^{\beta_2 } \} )BV( {\mathbb T}^2) = C( \{ n^{\beta_1 } \},\{ n^{\beta_2 } \} )V( {\mathbb T}^2)
\]
when $ \beta_1 = \beta_2 >\frac{1}{2}$.
Therefore, the following property is true.

\begin{consno}\label{cons-sum}
Let $m=2$, $\alpha=\alpha_1 = \alpha_2 \in (-\frac{1}{2},0)$ and $\beta = \alpha+1$. Then for any function $f$ from the class 
$ \{ n^{\beta } \} BV( {\mathbb T}^2)$ its Fourier series is $(C,\aaa)$ summable to $f^*( \xbf_0)$ at every
regular point $\xbf_0$.  If $f$ is continuous in a neighborhood of a compact set $K$, then summability is uniform on $K$.
\end{consno}

In section 4, we prove that for $ m \geq 3$ and certain $\{\alpha_j\}$ the continuity in variation 
is essential for the summability and even for the localization  of the Ces\'aro means (unlike the results on convergence,
see \cite{ownint} for details). 
More precisely, the following result is established.

\begin{thmno}\label{T:div}
Let $ m \geq 3$, $\alpha_j \in (-1,0)$ and  $\beta_j = \alpha_j+1$, $j=1,\dots,m$. 
Suppose that the condition
 \begin{equation}
\Bigl( \sum_{j=1}^m 	\beta_j \Bigr) - \beta_q >1
\label{bigbeta}
\end{equation}
 holds for a certain $q \in \{ 1,\dots,m\} $.
  Then there exists a continuous function $f$ with the following properties:
 
(1) $f$ belongs to the class $( \{ n^{\beta_1 } \},\dots, \{ n^{\beta_m } \} )BV( {\mathbb T}^m)$;

 (2) $f \equiv 0$ on $[-1,1]^m$;

 (3)  the cubic  $(C,\aaa)$ means of its Fourier series do not converge to zero at $\xbf_0=0$.
\end{thmno}

\begin{xrem}
The statement of Theorem \ref{T:div} holds for all $m\geq 2$ and $\alpha_j \in (-1,0)$, $j=1,\dots,m$ that
do not satisfy the conditions of Corollary \ref{cons-sum}.
But the proof for other $\{ \alpha_j \}$ is based on other methods. See \cite{owndan,ownzam} for details.
\end{xrem}

\section{Auxiliary estimates}\label{Saux}

First we recall some properties of the $\La$-variation obtained in our previous papers. 

\begin{lmmno}\label{onedim}
If $f$ is a function of $m$ variables, $f \in (\La^1,\dots,\La^m) BV(\Delta)$, 
$g$ is a function of one variable,
$g\in \La^j BV(\Delta^j)$, then
$F(\xbf) = f(\xbf) g(x^j) \in  (\La^1,\dots,\La^m) BV(\Delta)$
and the following estimate holds:
\begin{multline*}
V_{\La^1,\dots,\La^m} (F;\Delta) \leq C(m) \Bigl(
V_{\La^1,\dots,\La^m} (f;\Delta) V_{\La^j} (g; \Delta^j) +\\+
V_{\La^1,\dots,\La^m} (f;\Delta) \sup_{\Delta^j} | g | +
\sup_\Delta | f |\cdot  V_{\La^j} (g;\Delta^j)
\Bigr).
\end{multline*}
\end{lmmno}

This lemma was proved by the author in \cite[Lemma 2]{ownint} 
for $\La^1=\dots = \La^m$; in general case, the proof is almost the same.\qed

\begin{lmmno}\label{locvar}
Let a function $f(\xbf)$ belong to the class $C(\La^1, \dots,\La^m) V(\Delta)$.
Suppose there exists the limit $f(\xbf_0 + 0)$,
where $\xbf_0$ is an internal point of the interval $\Delta$. Then
\[
\lim_{\delta\to +0}
V_{ \La^1, \dots,\La^m }  \Bigl( f; \bigotimes_{j=1}^m (x^j_0, x^j_0 +\delta) \Bigr) =0.
\]
\end{lmmno}

This lemma was proved by the author in \cite[Theorem 2]{owncv}.\qed

In one-dimensional case, we introduce an auxiliary concept.

\begin{defn}
Let $\Delta \subset \mathbb T$ be an interval. 	For a set $E \subset \Delta$ put  by definition
\[
V_{\La} (f;\Delta\setminus E) = \sup_{\Omega (\Delta\setminus E)} \sum_k \frac{ |f(I_k) |}{\la_k}.
\]
\end{defn}

\begin{lmmno}
Let $\Delta \subset \mathbb T$  be an interval, $E \subset \Delta$ be a set with dense complement. 
Consider a function $f$ on $\Delta\setminus E$ such that  $V_{\La} (f;\Delta\setminus E)< \infty$. 
Then $f$ can be extended on $\Delta$ in such a way that
$V_{\La} (f;\Delta) = V_{\La} (f;\Delta\setminus E)$ and
$\sup_{\Delta} |f| = \sup_{\Delta\setminus E} |f|$.
\label{fullsegm}
\end{lmmno}

This lemma was proved by the author in \cite[Lemma 1]{owncoeft}.\qed

\begin{lmmno}\label{eqclass}
A function $f(\xbf)$ belongs to the class
$ C(\La^1,\dots,\La^m)V (\Delta)$,
if and only if there exist sequences $M^j=\{ \mu_n^j \} \in \LA$ such that 
$\dfrac{\mu^j_n}{\la^j_n} \downarrow 0$ as $n\to\infty$ and $f\in (M^1,\dots,M^m)BV(\Delta)$.
\end{lmmno}

This lemma was proved by the author in \cite[Theorem 1]{owncv}.\qed

We shall also use the following obvious properties of $\La$-variation.
For any functions $f$ and $g$ on an interval $\Delta$,
\begin{equation}
V_{\La^1,\dots, \La^m} (f+g; \Delta ) \leq
V_{\La^1,\dots, \La^m} (f; \Delta )  +
V_{\La^1,\dots, \La^m} (g; \Delta ).
\label{addfun}
\end{equation}
If $\Delta_1$ and $\Delta_2$ are two intervals,
their union is an interval and they intersect by a common face, then
\begin{equation}
V_{\La^1,\dots, \La^m} (f; \Delta_1 \cup \Delta_2 ) \leq
V_{\La^1,\dots, \La^m} (f; \Delta_1 )  +
V_{\La^1,\dots, \La^m} (f; \Delta_2 ).
\label{addset}
\end{equation}
It follows from \eqref{addset} that if $f$ is a $2\pi$-periodic function with respect to every variable and 
 $f \in (\La^1,\dots,\La^m) BV (\mathbb T^m)$, then for any $\xbf\in \mathbb R^m$ the following inequality holds:
\[
V_{\La^1,\dots, \La^m} (f; \xbf + \mathbb T^m ) \leq  2^m V_{\La^1,\dots, \La^m} (f;  \mathbb T^m ).
\]

Further on, we suppose that $\alpha_j \in (-1,0)$,
and $\beta_j = \alpha_j+1$, $j=1,\dots,m$.

\begin{lmmno}\label{lemB}
Let $m \geq 2$, $\{1,\dots,m\} = \{p \} \sqcup\, \xi$ for a certain $p$. Consider a function
$f\in (\La^1,\dots, \La^m) BV(\mathbb T^m)$,
where $\la^j_n = n^{\beta_j}$ for $j\in\xi$, and
$\Delta^p  \subset \mathbb T$, $\Delta^\xi \subset \mathbb T^{m-1}$ are intervals.
By definition, put 
\[
\varphi_\nn(t^p) =
\int_{\Delta^\xi}
f(x^p +t^p,x^\xi +t^\xi)
\prod_{j\in\,\xi} K^{\alpha_j}_{n_j} (t^j) \,dt^\xi.
\]
Here the integral exists for a.e. $t^p$, as
$f(x^p + t^p,x^\xi + \cdot )$ is bounded for all $t^p$ and is measurable for  a.e. $t^p$ by Fubini theorem.
Let $E = E (f,p, \xbf)$ be the set of Lebesgue measure zero where this integral does not exist.
Then for any $n_j>10$ the estimate
\[
V_{\La^p} (\varphi_{ \nn }; \Delta^p \setminus E ) \leq C(m,\aaa)
 V_{\La^1,\dots,\La^m} (f;  (x^p+\Delta^p)\times (x^\xi +\Delta^\xi) )
\]
holds.
\end{lmmno}

\begin{lmmno}\label{lemA}
Let $m \geq 1$, let $\Delta= \bigotimes\limits_{j=1}^m (a^j, b^j)$ be a subinterval in $\mathbb T^m$.
Then there exists a value $C(m,\aaa)$ such that
the estimate 
\begin{equation*}
\biggl|
\int_\Delta g(\tbf ) \prod_{j=1}^m K^{\alpha_j}_{n_j} (t^j) \,d\tbf
\biggr|
\leq C(m,\aaa)
\bigl( V_{ \{ n^{\beta_1} \}, \dots, \{ n^{\beta_m} \}  }
 ( g; \Delta ) +   \sup_{\Delta} |g|   \bigr)
\end{equation*}
holds for any  $g \in (\{ n^{\beta_1} \}, \dots, \{ n^{\beta_m} \} ) BV ( \Delta)$ and any $n_j>10$.
\end{lmmno}

\noindent{\itshape Proof of lemmas \ref{lemB} and \ref{lemA}.}
We prove these lemmas by induction on $m$.
First, we prove Lemma \ref{lemA} for $m=1$.

Let $\Delta= (a,b)$. We shall write $V_{\beta}$ instead of  $V_{ \{ n^\beta \} }$ for short.
We write $(t)^k$ for the $k$-th power of $t$ to avoid confusing it with the $k$-th coordinate
of a vector $\tbf$.
Without loss of generality we can assume that $a,b\geq 0$. 
The case $a,b \leq 0 $ is similar; if  $a<0$ and $b>0$, we represent the  integral
 as a sum of  integrals over $(a,0)$ and $(0,b)$ and estimate each of them.

Put $\eta_n =\frac{\pi}{n+ \bigl(\frac{\alpha+1}{2} \bigr) } $ for $ n >10$.
Our aim is to prove that
\begin{equation}
| A_n|  =\biggl|\, \int_a^b g(t) K_n (t)\, dt \biggr| \leq C(\alpha) ( \sup_{(a,\,b)} |g| + V_{\beta} (g;(a,b) )
\label{locbound}
\end{equation}
for any $a$ and $b$ such that $0\leq a <b \leq \pi$.

If $b \leq \eta_n$, then using the first of the estimates \eqref{estker},
we obtain
\begin{equation}
 | A_{n} | \leq 2n \frac{\pi}{n} \sup_{(a, b ) } |g| = 2\pi \sup_{(a, b ) } |g| . 
\label{estsmall}
\end{equation}
Let $b > \eta_n$. If $a < \eta_n$, then we split the integral into two:
\[
 A_n = \biggl(  \int_a^{\eta_n} +  \int_{\eta_n}^b \biggr)  g(t) K_n (t)\, dt =A_{n,1} + A_{n,2}.
\]
The first term can be estimated similarly to \eqref{estsmall}.
Thus, we should estimate the integral over $(a_n,b)$, where $a_n = \max \{ \eta_n ,a \}$.
In this integral, we decompose the Ces\'aro kernel  using \eqref{kernel}.
For the remainder term, we have:
\[
	\biggl|\,  \int_{a_n}^b g(t)
	\frac{ 2 \theta \alpha }{ n \bigl( 2 \sin \frac{t}{2} \bigr)^2  }\,dt \biggr| \leq
	\frac{C( \alpha) }{n}  \int_{a_n}^b  
	\frac{|g(t)| }{ (t)^2}\,dt 
	 \leq 
	\frac{C( \alpha)  }{n \eta_n } \sup_{ (a_n, b  ) } |g |
	 \leq C(\alpha) \sup_{ ( a_n, b  ) } |g |.
\]
For the main term, using the properties of sinus we obtain:
\[
I_n= \int_{a_n}^b g(t) K_n^{\alpha,*} (t)\, dt =  - \frac{1}{ A_n^\alpha } \int_{a_n - \eta_n}^{b-\eta_n} g(t+ \eta_n) 
\frac{ \sin \bigl[ \bigl( n + \frac{1}{2}+ \frac{\alpha}{2} \bigr) t
- \frac{\pi \alpha}{2}  \bigr] }{ \bigl( 2 \sin \frac{t+ \eta_n}{2} \bigr)^{\alpha+1 } }\,dt.
\]
Then
\begin{multline*}
2I_n = \int_{a_n}^{b-\eta_n} (g(t) - g(t+\eta_n)) K_n^{\alpha,*} (t)\, dt + \\+
\frac{1}{ A_n^\alpha } \int_{a_n}^{b-\eta_n} g(t+\eta_n) \sin [\dots] \biggl( \frac{1}{ (2 \sin \frac{t}{2})^{\alpha+1}}
- \frac{1}{ (2 \sin \frac{t+ \eta_n}{2})^{\alpha+1}} \biggr) \, dt + \\+
\int_{b-\eta_n}^b g(t) K_n^{\alpha,*} (t)\, dt +\\+
 \int_{a_n - \eta_n}^{a_n} g(t+\eta_n) K_n^{\alpha,*} (t+\eta_n)\, dt  = I_{n,1} + I_{n,2} + I_{n,3} + I_{n,4}.
\end{multline*}
We get using \eqref{estker} that
\[
|I_{n,3}| \leq 2\pi \sup_{(a,b)} |g|, \qquad |I_{n,4}| \leq 2\pi \sup_{(a,b)} |g|.
\]
On the other hand, for any  $ t \in (a_n, \pi)$ by the  mean value theorem we can find a point $\tau \in (t, t+ \eta_n )$ 
such that 
\[
\biggl|
\frac{1}{ (2 \sin \frac{t}{2})^{\alpha+1}}
- \frac{1}{ (2 \sin \frac{t+ \eta_n}{2})^{\alpha+1}} \biggr| =
\biggl| - \eta_n \frac{  (\alpha +1) \cos \frac{\tau}{2} }{ (2 \sin \frac{\tau}{2} )^{\alpha+2}}  \biggr|
\leq
\frac{ C \eta_n }{ (t)^{\alpha+2} }.
\]
As $a_n \geq \eta_n$, the second term can be estimated in the following way:
\[
|I_{n,2}| \leq \frac{C   }{n^\alpha}  \eta_n \int_{a_n}^{b-\eta_n} \frac{|g(t+\eta_n)|\, dt}{ (t)^{\alpha+2}} 
\leq 
\frac{C}{n^\alpha} \frac{\eta_n}{ (a_n)^{\alpha+1} } \sup_{ (a_n,b ) } |g| \leq 
C \sup_{ (a_n,b) } |g| .
\]
Finally, for the first term we have:
\[
|I_{n,1}| \leq 
\frac{C}{ A_n^\alpha } 
\int_{a_n}^{b-\eta_n} \frac{ |g(t) - g(t+\eta_n) |}{ (t)^{\alpha+1 } }\, dt = J_n.
\]
Let $k$ be the integral part of $\frac{b- a_n}{\eta_n}-1$. Put $\varkappa_n = (b-\eta_n) -  a_ n - k \eta_n$. 
Then we have
\begin{multline*}
	J_n = \frac{C}{ A_n^\alpha } \biggl( \sum_{j=1}^{k}  \int_{a_n +(j-1)\eta_n}^{a_n +j\eta_n} \frac{ |g(t) - g(t+\eta_n) |}{ (t)^{\alpha+1 } }\, dt
	 \\ + 
	 \int_{a_n+ k \eta_n}^{a_n + k \eta_n + \varkappa_n} \frac{ |g(t) - g(t+\eta_n) |}{ (t)^{\alpha+1 } }\, dt \biggr).
\end{multline*}
Taking into account that  $a_n + (j-1) \eta_n \geq j \eta_n$, we get
\begin{multline*}
	J_n \leq \frac{C}{ n^\alpha } 
	\biggl( \int_{0}^{\varkappa_n} \biggl( \sum_{j=1}^{k+1} \frac{ |g(a_n+ (j-1) \eta_n+t) - g(a_n+j\eta_n+t) |}{ (j \eta_n)^{\alpha+1 } } \biggr)\, dt 
 \\ + 
	\int_{\varkappa_n}^{\eta_n} \biggl( \sum_{j=1}^{k} \frac{ |g(a_n+ (j-1) \eta_n+t) - g(a_n+j\eta_n+t) |}{ (j \eta_n)^{\alpha+1 } } \biggr)\, dt 
	\biggr).
\end{multline*}
For every $t$, both integrand expressions equal a variation sum multiplied by $(\eta_n)^{-(\alpha+1)}$.
Taking into account that $\eta_n \sim  \frac{1}{n}$, we get
\[
J_n \leq \frac {C \eta_n V_{\beta   } (g;(a,b)) }{ n^\alpha (\eta_n)^{\alpha+1} }  \leq C V_{\beta   } (g; (a,b) ).
\]
Therefore, the estimate \eqref{locbound} holds, and 
Lemma \ref{lemA} is proved for $m=1$.

Now consider that  lemma \ref{lemA} is proved for $(m-1)$-dimensional case, $m \geq 2$.
Let us prove Lemma \ref{lemB} for $m$-dimensional case.
Without loss of generality we can assume that $p=1$.
Let $F(\tbf) = f(\xbf+\tbf)$.
Consider a system of intervals $\{ I^1_{k_1} \}\in\Omega(\Delta^1 \setminus E)$.
For any $ z \in \mathbb C$, denote 
$\sgn z =  |z|/z$, $z \neq 0$, and $\sgn 0 =0$.
Putting $\eps_{k_1} =
\sgn( \varphi_{\nn} ( I^1_{k_1}  ))$, we have
\begin{equation}
Z =
\sum_{ k_1 }
\frac{|\varphi_{\nn} ( I^1_{k_1} ) | }{ \la^1_{k_1}  } =
\sum_{k_1 } \frac{\eps_{k_1}
\varphi_{\nn} ( I^1_{k_1} )  }{ \la^1_{k_1} }
=
\int_{\Delta^\xi} \psi(t^\xi)
\prod_{j\in\xi} K_{n_j}^{\alpha^j} (t^j) \, dt^\xi,
\label{varsum}
\end{equation}
where the function $\psi$ is defined by the formula
\begin{equation*}
\psi(t^\xi) =
\sum_{ k_1 }
\frac{\eps_{k_1} }{ \la^1_{k_1} } F(I^1_{k_1} , t^\xi).
\end{equation*}
Using the definition of $\La$-variation, we have
\[
\sup_{\Delta^\xi} |\psi| \leq V_{\La^1}^{x^1} (F; \Delta ) \leq V_{\La^1,\dots,\La^m} (F; \Delta ).
\]
Our task is now to estimate the $(\La^2 ,\dots,\La^m)$-variation of the function $\psi$.
Consider a nonempty subset $\tau \subset \xi$, and let $\varkappa = \xi\setminus \tau$.
For any systems of intervals
$\{ I^j_{k^j} \} \in\Omega(\Delta^j)$, $j\in\tau$ we have
\[
\sum_{ k^\tau }
\frac{| \psi(I^{\tau}_{k^{\tau}} ) |}%
{ \la^\tau_{k^\tau} } \leq
\sum_{k_1, k^\tau}
\frac{ |F (I^1_{k_1} \times I^\tau_{k^\tau},x^\varkappa )| }%
{\la^1_{k_1} \la^\tau_{k^\tau} }
\leq
V_{\La^1, \La^\tau}^{ x_1, x^\tau}  (F; \Delta).
\]
Taking the supremum over $\Omega (\Delta^j) $, $j\in\tau$
and $x^\varkappa \in \Delta^\varkappa$, and then summing over $\tau$,
we obtain
\[
V_{\La^\xi}(\psi;\Delta^\xi) \leq  V_{\La^1,\dots,\La^m} (F;\Delta).
\]
Applying Lemma \ref{lemA} to the function $\psi$ (of $(m-1)$ variables) 
and taking into account (\ref{varsum}), we get
\[
Z \leq C(m-1,\aaa) (V_{\La^\xi} (\psi; \Delta^\xi) +\sup_{\Delta^\xi} |\psi | ).
\]
Taking the supremum over $\Omega (\Delta^1 \setminus E) $, 
we obtain 
\[
V_{\La^1} (\varphi_\nn; \Delta^1 \setminus E)
\leq  C(m,\aaa) V_{\La^1,\dots,\La^m}(F; \Delta)
= C(m,\aaa) V_{\La^1,\dots,\La^m} (f; \xbf+\Delta ).
\]
Thus Lemma  \ref{lemB} is proved for $m$-dimensional case.

Now suppose that $m \geq 2$, Lemma \ref{lemB} is proved  for $m$-dimensional case
and Lemma \ref{lemA} is proved for $(m-1)$-dimensional case.
Let us prove Lemma \ref{lemA} for $m$-dimensional case.
Put $ t^* =(t^2,\dots,t^m)$ and
$\Delta^*= \bigotimes\limits_{j=2}^m (a^j, b^j)$.
We have
\[
	Z_1=
	\int_\Delta g(\tbf) \prod_{j=1}^m K_{n_j}^{\alpha_j} (t^j) \,d\tbf
	=
	\int_{a^1}^{b^1} K_{n_1}^{\alpha_1} (t^1) 
	\biggl(\,
	\int_{ \Delta^*} g(t^1, t^*) \prod_{j=2}^m K_{n_j}^{\alpha_j} (t^j)  \,dt^*
	\biggr)\,dt^1.
\]
Consider the function
\[
G(t^1)=
\int_{\Delta^*} g(t^1,  t^*)
\prod_{j=2}^m K_{n_j}^{\alpha_j} (t^j)  \,d{ t^* }.
\]
Let $\Lambda^j = \{ n^{\beta_j} \}$ for $j=1,\dots,m$.
By Lemma \ref{lemB} we have
$
G(t^1)\in \La^1 BV( \Delta^1 \setminus E)$  and 
$V_{\La^1} (G; \Delta^1 \setminus E) \leq C_1(m,\aaa)  V_{\La^1,\dots,\La^m} ( g;\Delta )
$,
where a set $E$ is of Lebesgue measure zero. Hence $\Delta^1 \setminus E$ is dense on $\Delta^1$.
By the inductive hypothesis,
\[
\sup_{\Delta^1\setminus E} | G |
\leq C(m-1,\aaa) \bigl( V_{\La^1,\dots,\La^m} (g; \Delta)+ \sup_{\Delta} |g| \bigr) .
\]
By Lemma \ref{fullsegm} the function $G$ can be extended on $\Delta^1$ without increasing of its $\La^1$-variation
and supremum.
Then, applying Lemma \ref{lemA} for $m=1$, we get
\[
|Z_1| \leq C(1,\aaa)
\bigl(V_{\La^1} (G;\Delta^1) + \sup_{\Delta^1} |G | \bigr).
\]
Combining these estimates, we obtain
\[
|Z_1| \leq C(m,\aaa) \bigl( V_{\La^1,\dots,\La^m} (g; \Delta) +
  \sup_{\Delta} |g | \bigr).
\]
Lemmas \ref{lemB} and \ref{lemA} are proved.\qed

\begin{lmmno}\label{finint}
Let $ \Delta =[a,b] \subset \mathbb R$ and
 $f \in \La BV ( \Delta)$. Let a function $s (t)$ be continuous and satisfy  the conditions $|s(t)| \leq 1$,
 $s (t+\pi) = - s (t)$. 
   Then for $\varkappa = ( b - a)/ \pi $ the estimate 
\[
\biggl| \int_\Delta f(t) s (A t) \,d t\biggr|
\leq
 C( \Delta )   \biggl( \frac{ V_{\La} (f;\Delta) }%
{ \La( \lfloor \varkappa A \rfloor -1 ) }
+
  \frac { \sup_\Delta |f| }{A}  \biggr).
\]
holds  for any $A \geq \frac{2}{\varkappa}$ (hereafter by $\lfloor x \rfloor $ we denote the integral part of $x$).
\end{lmmno}

The lemma was proved in \cite[Lemma 6]{ownint} for  $s(t) = \sin t$.
In the general case, the proof is almost the same.

\section{Proof of the result on summability} \label{Sconverge}

\noindent{\itshape Proof of Theorem \ref{T:main}.}
Let $\alpha_j \in (-1,0)$ and $\beta_j = \alpha_j+1$, $j=1,\dots,m$.
Consider a function $f \in C( \{ n^{\beta_1 } \},\dots, \{ n^{\beta_m } \} )V (\mathbb T^m)$.
Applying Lemma \ref{lemB} for $\Delta = \mathbb T^m$,
we obtain the first statement of the theorem (the boundedness of Ces\'aro means).

Now we shall prove the second statement. Consider the function
\[
\varphi_{\xbf} (\tbf) = \frac{1}{2^m} \sum f(x^j \pm t^j) - f^* (\xbf).
\]
From \eqref{addfun} we get that this function belongs to any Waterman class that $f$ does, and
\begin{equation}
 V_{\La^1,\dots,\La^m} (\varphi_\xbf; [0,\pi]^m ) \leq C(m)  V_{\La^1,\dots,\La^m} (f; \mathbb T^m)
\label{eq:phix}
\end{equation} 
uniformly with respect to $\xbf$.
The Ces\'aro means of $f$ have the form 
\[
\sigma_{\nn}^{\aaa} (\xbf) - f^* (\xbf) = \frac{2^m}{\pi^m} \int_{ [0,\pi]^m} \varphi_{\xbf} (\tbf) 
\prod_{j=1}^m K_{n_j}^{\alpha_j} (t^j) \,d\tbf.
\]
Take an $\eps>0$.
First, by Lemma \ref{locvar} we can take a $\delta>0$ such that
\[
V_{ \{ n^{\beta_1 } \},\dots, \{ n^{\beta_m } \} } (\varphi_\xbf;  (0, \delta)^m )
\leq \frac{\pi^m \eps}{2 C (m,\aaa)},
\]
where $C(m,\aaa)$ is from Lemma \ref{lemA}. 
If the function $f$ is continuous, $\delta$ does not depend on a point $\xbf$.
For an arbitrary partition of the set $\{1,\dots,m \}$
into two non-intersecting subsets $\gamma$ and $\xi$
we put
\[
P_{\gamma,\xi}
= \bigotimes_{j=1}^m J^j_{\gamma,\xi},
\mbox{ where }
J^j_{\gamma,\xi} =
\begin{cases}
( 0, \delta ), & j\in\gamma , \cr
[ \delta, \pi ), & j\in\xi .
\end{cases}
\]
Then we get
\begin{equation}
\tfrac{\sigma_{\nn}^{\aaa} (\varphi_\xbf,\zero) - f^* (\xbf)}{2^m}  =
\sum_{\gamma \sqcup \xi = \{ 1,\dots,m \} }
S_{\nn}^{\gamma,\xi } =
\sum_{\gamma \sqcup \xi = \{ 1,\dots,m \} }
\frac{1}{\pi^m}
\int_{ P_{\gamma,\xi} } \varphi_\xbf (\tbf) \prod_{j=1}^m K_{n_j}^{\alpha_j} (t^j) \,d\tbf.
\label{thesum}
\end{equation}

If $\xi = \varnothing$, then using Lemma \ref{lemA} we obtain the inequality 
$| S_\nn^{\gamma, \varnothing}| < \eps/2 $.
We shall show that all other terms in \eqref{thesum} tend to zero as $\nn$ grows. 
Consider a nonempty $\xi$. Without loss of generality,
$\xi = \{ p+1, \dots, m \}$ for a certain $p<m$.
Put $\La^j = \{ n^{\beta_j } \}$, $j=1,\dots,m-1$.
Using Lemma \ref{eqclass}, we find a sequence $\La^m$ such that
$\frac{\la^m_k}{k^{\beta_m}} \to 0$ as $k\to\infty$ and
$f\in (\La^1,\dots, \La^m) BV (\mathbb T^m)$.
Let
\[
F_{\nn} (t^m) = \int_{[0,\delta]^p} \int_{(\delta,\pi]^{m-p-1}} \varphi_\xbf (\tbf) \prod_{j=1}^{m-1}
K_{n_j}^{\alpha_j} (t^j) \, d( t^1 \dots  t^{m-1})
\]
In view of \eqref{eq:phix} by Lemma \ref{lemA} we have 
\[
|F_{\nn} (t^m)| \leq C(m-1,\aaa) \Bigl( V_{\La^1,\dots, \La^m} (f; P_{\gamma,\xi}) + \sup_{ P_{\gamma,\xi} } |f| \Bigr)
\]
uniformly with respect to $\xbf$ and $\nn$.
By Lemma \ref{lemB}, these functions belong to the class $\La^m BV ([\delta,\pi])$, and
using \eqref{eq:phix} we obtain
 \[
V_{\La^m} (F_{\nn};  [\delta,\pi]) \leq C(m,\aaa)  V_{\La^1,\dots, \La^m} (f; P_{\gamma,\xi})
\]
uniformly with respect to $\xbf$ and $\nn$.
Consider the functions 
\[
G_\nn(t^m) = F_\nn (t^m) \cdot  \frac{\chi_{[\delta,\pi]} (t^m) }{ (2\sin\frac{t^m}{2})^{\beta_m} }.
\]
By Lemma \ref{onedim} with regard to \eqref{eq:phix} these functions belong to the class $\La^m BV (\mathbb T)$,
and the estimates
\begin{gather}
	\sup_{t\in \mathbb T} | G_\nn (t) | < C (m,\aaa,\delta) \Bigl( V_{\La^1,\dots, \La^m} (f; \mathbb T^m) + \sup_{ \mathbb T^m } |f| \Bigr),
	\label{g1}\\
 V_{\La^m} (G_{\nn};  [\delta,\pi]) \leq C (m,\aaa,\delta) \Bigl( V_{\La^1,\dots, \La^m} (f; \mathbb T^m) + \sup_{ \mathbb T^m } |f| \Bigr).	
 \label{g2}
\end{gather}
hold uniformly with respect to  $\nn$ and $\xbf$.
Put $\nu_n = n+ \frac{\alpha_m+1}{2} $.
Then
\begin{multline}
\pi^m S_\nn^{\gamma,\xi } = \int_\delta^\pi F_\nn (t^m) K_{n^m}^{\alpha_m} (t^m) \, dt^m
=  \frac{ \cos \frac{\pi\alpha_m}{2}  }{A_{n_m}^{\alpha_m} }  \int_{0}^\pi G_\nn (t^m) \sin (\nu_{n_m} t^m) \,dt^m - \\ -
 \frac{ \sin \frac{\pi\alpha_m}{2}  }{A_{n_m}^{\alpha_m} }  \int_{0}^\pi G_\nn (t^m) \cos (\nu_{n_m} t^m) \,dt^m +
 \int_\delta^\pi F_\nn (t^m) R_{n^m}^{\alpha_m} (t^m) \, dt^m.
\end{multline}
Here the last term can be estimated as follows:
\[
\biggl|\, \int_\delta^\pi F_\nn (t^m) R_{n^m}^{\alpha_m} (t^m) \, dt^m \biggr| \leq \sup_{ [\delta,\pi] } | F_\nn (t) | \cdot \frac{ 2|\alpha_m| }{n_m}
\cdot \int_\delta^\pi \frac{ dt}{ (2 \sin \frac{t}{2} )^2 }\to 0 
\]
as $n_m \to\infty$, independently of other components of $\nn$.

Consider the first term (the second is estimated likewise).
By Lemma~\ref{finint},
\[
\biggl|\, \int_{0}^\pi G_\nn (t^m) \sin (\nu_{n_m} t^m) \,dt^m \biggr|
\leq C \biggl(
\frac{ V_{\La^m} (G_\nn; [\delta,\pi]) }{ \La^m ( \lfloor \nu_{n_m} \rfloor -1 )}
  + \frac{ \sup_{t\in \mathbb T} | G_\nn (t) | }{ \nu_{n_m} }
\biggr).
\]
Recalling \eqref{g1} and \eqref{g2}, we get
\begin{multline*}
	\frac{ 1 }{A_{n_m}^{\alpha_m} } \biggl|\, \int_{0}^\pi G_\nn (t^m) \sin (\nu_{n_m} t^m) \,dt^m \biggr| 
	\\ \leq 
	C( m,f, \aaa, \delta ) \biggl( \frac{1}{ A_{n_m}^{\alpha_m}  \La^m ( {n_m} -1 )} + \frac{1 }{ A_{n_m}^{\alpha_m} \nu_{n_m} } \biggr).
\end{multline*}
Here $A_{n_m}^{\alpha_m} \sim (n_m)^{\alpha_m}$, $ \nu_{n_m} \sim n_m$, hence $A_{n_m}^{\alpha_m} \nu_{n_m} 
\sim (n_m)^{\beta_m} \to \infty$ as $n_m \to\infty$.
Finally, we have chosen  $\La^m$ such that the condition
$\frac{\la^m_k}{k^{\beta_m}} \to 0$ as $k\to\infty$
holds, therefore,
$
\frac{ \La^m (N) }{ \sum_{k=1}^N \frac{1}{ k^{\beta_m} } }\to \infty,
$
where $\sum_{k=1}^N \frac{1}{ k^{\beta_m} } \sim \frac{1}{N^{\alpha_m}}$. Consequently, 
$  A_{n_m}^{\alpha_m}  \La^m ( n_m -1 ) \to \infty$
as $n_m \to\infty$.
This completes the proof of Theorem \ref{T:main}.
\qed

\section{Proof of the result on non-summability}\label{S:diverg}

We use the following construction introduced and studied in our paper \cite{owncv}.
Let $m\geq 3$.
Consider a system $\{ \JJ^1_k \}_{k=1}^\infty$ of intervals enclosed into $\mathbb T$
but not covering $\mathbb T$.
Consider systems  $\{ \JJ^j_k \}_{k=1}^\infty$ of non-intersecting intervals,
$j=2,\dots,m$, where $\JJ^j_k = (a^j_k,b^j_k)\subset \mathbb T$.

Let $f_k(x)$ be arbitrary functions on
$\mathbb T$ such that $f_k(t)=0$ for
$t\leq a^1_k$ and  $f_k(t)=0$ for $t\geq b^1_k$.
Let  $h_j^k (t)$ be arbitrary functions on
$\mathbb T$ such that 
$h_j^k (t)=0 $ for $t\leq a_k^j$ and $t\geq b_k^j$, 
 $h_k^j (( a_k^j+b_k^j )/2)= 1$, $h_k^j(t)$ do not decrease on
$ [a_k^j;( a_k^j + b_k^j )/2]$
and they do not increase on $[ ( a_k^j + b_k^j )/2; b_k^j ]$.

We say that  $f$ is a ``diagonal'' function on $\mathbb T^m$ if it is the sum of the series
\begin{equation}
f(\xbf) = \sum_{k=1}^\infty
\Bigl( f_k (x^1 )  \prod_{j=2}^m h_k^j(x^j) \Bigr)
\label{diag}
\end{equation}
where $f_k$ and $h_k^j$ are described above.
It is obvious that the support of such a function is contained 
in the union of closures of the pairwise non-intersecting parallelepipeds
$\JJ^1_k \times \dots \times \JJ^m_k$.

\begin{lmmno}[{\cite[Lemma 6]{owncv}}]\label{diagf}
Consider sequences $\La^2,\dots \La^m \in \LA$ such that
\begin{equation}
\sum_{k=1}^\infty \frac{1}{\la^2_k \dots \la^m_k} <\infty.
\label{ladiv}
\end{equation}
Then the following conditions are equivalent:

(a) the ``diagonal'' function $f(\xbf)$ defined by the formula (\ref{diag})
belongs to the class  $(\La^1,\dots,\La^m)BV( \mathbb T ^m)$;

(b) $f_k\in \La^1 BV( \mathbb T)$ for every $k\in\mathbb N$,
and  $\sup_k V_{\La^1}(f_k; \mathbb T ) =C <\infty$.
\end{lmmno}

\noindent
{\itshape Proof of Theorem \ref{T:div}.}
Suppose $\alpha_j \in (-1,0)$ and
$\beta_j = \alpha_j+1$, $j=1,\dots,m$ satisfy \eqref{bigbeta}. 
Without loss of generality we may assume that $q=1$.
By definition, put 
\[
\rho = \frac{1}{2\pi \cdot 4^m} \int_1^3 \frac{dt}{t^{\beta_1  }}.
\]
We shall define inductively an increasing sequence of positive integers $\{ N_k \}$
and $3(m-1)$ sequences $c^j_k$, $d^j_k$, and $\delta^j_k$ ($j=2,\dots,m$).
Let $\nu_{j,k} = N_k + \frac{1+\alpha_j}{2}$.
Let $a_k $ and  $b_k$ be the smallest and the largest zeros of the function $\sin (\nu_{1,k} t - \frac{\pi \alpha_1}{2}) $  on the segment $[1,3]$.
Let
\[
f_k (t) = \chi_{[a_k,b_k]} (t) \cdot \sin \biggl(\nu_{1,k} t - \frac{\pi \alpha_1}{2} \biggr) \cdot  A_{N_k}^{\alpha_1} 
\]
We take $N_1=10$, $c^j_1 = 1/4$, $d^j_1 =1/2$, $\delta^j_1 = 1/2$.
Define  $h_1^j(t)$ as follows: $h_1^j=0$ outside $(1/4,1/2)$,  $h_1^j(3/8)=1$ and it is
linear on the segments $[1/4,3/8]$ and $[3/8,1/2]$.

Suppose that $N_k$, $c^j_k$, $d^j_k$ and $\delta^j_k$ are already defined for $k<s$. 
Let
\[
\psi_k (\xbf) = f_k (x^1 )  \prod_{j=2}^m h_k^j(x^j),\qquad
g_s (\xbf) = \sum_{k=1}^{s-1} \psi_k(\xbf) .
\]
It can easily be checked that the function $g_s$ satisfies the conditions of Theorem \ref{T:main}
for any functions $h^j_k \in BV ( \mathbb T)$.
Hence, $\sigma_{\nn}^{\aaa} ( g_s,\zero ) \to 0$ as $ \nn \to \infty$. 
Therefore, there exists  $M_{s,1}$ such that for any $n>M_{s,1}$ we have
$	| \sigma_{n}^{\aaa} (g_s,\zero)  | < \rho/4$.  
Let 
\[
\delta^j_s = \min \Bigl\{ \frac{1}{2} \delta^j_{s-1}, c^j_{s-1} , \frac{1}{ 4 N_{s-1}} \Bigr\}.
\]
Using the equality \eqref{unit} and the second of the estimates \eqref{estker} we find $M_{s,j}$, $j=2,\dots,m$ such that
\[
\frac{1}{\pi} \int_0^{\delta^j_{s}} K_n^{\alpha_j}(t)\, dt > \frac{5}{12}
\]
for every $j$ and for all $n> M_{s,j}$.
Let us take $N_s > \max_{j=1,\dots,m} M_{s,j}$ such that
\begin{equation}
	 (N_s/N_{s-1})^{\alpha_1} <  \min \Bigl\{ \frac{1}{8 B(\alpha_1)}, \frac{1}{2} \Bigr\},
\label{ngrowth}
\end{equation}
where $B(\alpha_1)$ is the constant in the second of inequalities \eqref{estker}.
Then using the properties of Lebesgue integral, we can take points $c^j_s$ and $d^j_s$ ($j=2,\dots,m$) 
such that $[c^j_s, d^j_s] \subset (0, \delta^j_s)$ and
\[
\frac{1}{\pi} \int_{c^j_{s}}^{d^j_{s}} K_{N_s}^{\alpha_j}(t)\, dt > \frac{1}{3}.
\]
After it, we take continuous functions $h^j_s$ with support on $[c^j_s, d^j_s]$ 
satisfying the conditions of Lemma \ref{diagf} (see the third paragraph of section~\ref{S:diverg}) 
and, therefore, belonging to $BV (\mathbb T)$ such that
\begin{equation}
	\sigma_{N_k}^{\alpha_j} (h^j_k,0) = \frac{1}{\pi} \int_{c^j_{s}}^{d^j_{s}} h^j_s(t)  K_{N_s}^{\alpha_j}(t)\, dt > \frac{1}{4}.
	\label{bigh}	
\end{equation}
It can be easily seen that 
\[
V_{ \{ n^{\beta_1} \} } (f_k; [0,\pi]) \leq A_{N_k}^{\alpha_1} 
\sum_{l=1}^{\lfloor 2N_k/\pi \rfloor} \frac{1}{l^{\beta_1}} \leq C (\alpha_1).
\]
Using the condition \eqref{bigbeta}, we see that the class 
$( \{ n^{\beta_1 } \},\dots , \{ n^{\beta_m } \} )BV( \mathbb T^m)$
satisfies \eqref{ladiv}.
Therefore, the ``diagonal'' function $f$, defined by \eqref{diag} satisfies the conditions of Lemma \ref{diagf}
and, consequently, belongs to the considered Waterman class. As each function $\psi_k(\xbf)$ is continuous, 
 their supports are pairwise non-intersecting and $\sup |\psi_k| \to 0$ as $k \to \infty$,
 the series \eqref{diag} converges uniformly and thus $f$ is continuous.

We now estimate the Ces\'aro means of its Fourier series.
First, for any $s$ we have:
\begin{multline*}
	\sigma_{N_s}^{\alpha_1} (f_s,0) = \frac{ A_{N_s}^{\alpha_1} }{\pi} \int_{a_s}^{b_s} \sin \biggl(\nu_{1,s} t - \frac{\pi \alpha_1}{2} \biggr)
	K_{N_s}^{\alpha_1} (t)\, dt
 	=
	\frac{ 1 }{\pi} \int_{a_s}^{b_s} \frac{ \sin^2 \bigl(\nu_{1,s} t - \frac{\pi \alpha_1}{2} \bigr) }{(2\sin t/2)^{\alpha+1} }\, dt
	\\+
	  \frac{ A_{N_s}^{\alpha_1} }{\pi} \int_{a_s}^{b_s} \sin \biggl(\nu_{1,s} t - \frac{\pi \alpha_1}{2} \biggr)
	R_{N_s}^{\alpha_1} (t)\, dt = J_1+ J_2.
\end{multline*}
The first term can be estimated in this way:
\begin{multline*}
	J_1 > \frac{1}{\pi} \int_{a_s}^{b_s} \frac{ \sin^2 \bigl(\nu_{1,s} t - \frac{\pi \alpha_1}{2} \bigr) }{t^{\alpha+1} }\, dt
	\\=
	 \frac{1}{2\pi} \int_{a_s}^{b_s} \frac{dt}{t^{\alpha+1} }  - 
	\frac{1}{2\pi} \int_{a_s}^{b_s} \frac{ \cos \bigl( 2\nu_{1,s} t - {\pi \alpha_1} \bigr) }{t^{\alpha+1} }\, dt =
J_{1,1} - J_{1,2}.
\end{multline*}
Here $J_{1,1}$ tends to $4^m \rho$ as $s$ tends to infinity, and $J_{1,2}$ tends to zero 
as $s$ tends to infinity by the mean value theorem.
On the other hand, 
\[
|J_2| \leq \frac{  A_{N_s}^{\alpha_1} }{ N_s \pi} \int_1^3 \frac{dt}{ (2\sin t/2)^2 } \to 0 
\]
as $s \to\infty$. Therefore, if $s$ is sufficiently large then
$\sigma_{N_s}^{\alpha_1} (f_s, 0)  > 4^m \rho /2$.
Combining this with \eqref{bigh}, we obtain that
\[
|\,\sigma_{N_s}^{\aaa} (\psi_s,\zero) | > \frac {4^m \rho }{2\cdot 4^{m-1}} > \rho
\]
for sufficiently large $s$.

Secondly, consider that $k>s$. Then the second of the estimates \eqref{estker} and the condition \eqref{ngrowth} imply the inequality
\[
	|\sigma_{N_s}^{\alpha^1} (f_k,0)| \leq (N_k)^{\alpha^1} \int_1^3 \frac{ B(\alpha^1) }{ (N_s)^{\alpha^1} t^{ \beta^1} }\, dt
	\leq  \frac{\rho}{8}.
\]
Moreover, for $j=2,\dots,m$, using the first of the estimates \eqref{estker} we obtain
\[
	|\,\sigma_{N_s}^{\alpha^j} (h^j_k,0)|  =\biggl| \frac{1}{\pi} \int_{c^j_{k}}^{d^j_{k}} K_{N_s}^{\alpha^j}(t)\, dt  
	  \biggr| \leq 2 N_s (d^j_k - c^j_k) \leq \frac{2}{\pi} N_s \delta^j_k.
\]
Combining these estimates and taking into account the conditions posed on $\delta^j_k$,
 we come for any $k>s$ to the inequality
\[
|\,\sigma_{N_s}^{\aaa} (\psi_k,\zero) | \leq \frac{ \rho }{8} \prod_{j=2}^m \biggl(  N_s \delta^j_{s+1}  \frac{ \delta^j_k }{ \delta^j_{s+1} } \biggr)
\leq \frac{\rho}{8} \biggl( \frac{1}{2} \biggr)^{(m-1)(k-s-1)}.
\]
Finally, as the series 
\[
f(\xbf) = \sum_{k=1}^\infty \psi_k (\xbf)
\]
converges uniformly on $\mathbb T^m$, we obtain that
\begin{multline*}
	|\,\sigma_{N_s}^{\aaa} (f,\zero) | \geq |\,\sigma_{N_s}^{\aaa} (\psi_s,\zero) | - |\,\sigma_{N_s}^{\aaa} (g_s,\zero) | - \sum_{k=s+1}^\infty
	|\,\sigma_{N_s}^{\aaa} (\psi_k,\zero) |
	\geq \\ \geq
	\rho - \frac{\rho}{4} - \frac{\rho}{8} \sum_{k=s+1}^\infty \biggl( \frac{1}{2} \biggr)^{(m-1)(k-s-1)} \geq \frac{\rho}{2}
\end{multline*}
for sufficiently large $s$.
This completes the proof of Theorem \ref{T:div}.\qed

{\small
{\em Author's addresses}:
{\em Alexandr Bakhvalov}, Department of Mechanics and Mathematics, Moscow Lomonosov University, Moscow, Russia; 
 e-mail: \texttt{an-bakh@\allowbreak yandex.ru}.

}

\end{document}